\documentclass[11pt,dvips]{article}
\usepackage{amsmath,amsfonts,amssymb,amsthm,graphicx,verbatim,subfig}
\usepackage[all]{xy}
\usepackage[a4paper,left=25mm]{geometry}
\font\sixbb=msbm6
\font\eightbb=msbm8
\font\twelvebb=msbm10 scaled 1095
\newfam\bbfam
\textfont\bbfam=\twelvebb \scriptfont\bbfam=\eightbb
                           \scriptscriptfont\bbfam=\sixbb
\def\bb{\fam\bbfam\twelvebb}
\newcommand{\Rea}{{\bb R}}

\newcommand{\Int}{{\bb Z}}
\newcommand{\Rat}{{\bb Q}}

\newtheorem{theorem}{\bf Theorem}
\newtheorem{claim}[theorem]{\bf Claim}
\newtheorem{conjecture}[theorem]{\bf Conjecture}

\newtheorem{corollary}[theorem]{\bf Corollary}
\newcommand{\enp}{\begin{flushright} $\Box$ \end{flushright}}
\newcommand{\beq}[0]{\begin{equation}}
\newcommand{\enq}[0]{\end{equation}}

\newcommand{\lk}{{\rm lk}}
\newcommand{\st}{{\rm st}}

\title{A Tverberg Type Theorem for Matroids}

\author{Imre B{\'a}r{\'a}ny \and  Gil Kalai \and Roy Meshulam}
\begin{document}
\maketitle
\begin{abstract}
Let $b(M)$ denote the maximal number of disjoint bases in a matroid
$M$. It is shown that if $M$ is a matroid of rank $d+1$, then for
any continuous map $f$ from the matroidal complex $M$ into $\Rea^d$
there exist $t \geq \sqrt{b(M)}/4$ disjoint independent sets
$\sigma_1,\ldots,\sigma_t \in M$ such that $\bigcap_{i=1}^t
f(\sigma_i) \neq \emptyset$.
\end{abstract}

\section{Introduction}

\ \ \ \ Tverberg's theorem \cite{Tv} asserts that if $V \subset \Rea^d$
satisfies $|V| \geq (k-1)(d+1)+1$, then there exists a partition
$V=V_1 \cup \cdots \cup V_k$ such that $\bigcap_{i=1}^k
\text{conv}(V_i) \neq \emptyset$. Tverberg's theorem and some of its
extensions may be viewed in the following general context. For a
simplicial complex $X$ and $d \geq 1$, let the {\it affine Tverberg
number} $T(X,d)$ be the maximal $t$ such that for any affine 
map $f: X \rightarrow \Rea^d$, there exist disjoint simplices
$\sigma_1,\ldots,\sigma_t \in X$ such that $\bigcap_{i=1}^t f(\sigma_i)
\neq \emptyset$. The {\it topological Tverberg number} $TT(X,d)$ is
defined similarly where now $f:X \rightarrow \Rea^d$ can be an
arbitrary continuous map.

Let $\Delta_n$ denote the $n$-simplex and let $\Delta_n^{(d)}$ be
its $d$-skeleton. Using the above terminology, Tverberg's theorem is
equivalent to $T(\Delta_{(k-1)(d+1)},d)=k$ which is clearly the same as $T(\Delta_{(k-1)(d+1)}^{(d)},d)=k$. Similarly, the topological Tverberg theorem of B{\'a}r{\'a}ny, Shlosman and Sz\H{u}cs \cite{BSS81} states that if $p$ is prime then $TT(\Delta_{(p-1)(d+1)},d)=p$. Sch\"{o}neborn and Ziegler \cite{SZ05} proved that this implies the stronger statement $TT(\Delta_{(p-1)(d+1)}^{(d)},d)=p$. This result was extended by \"Ozaydin~\cite{Oza} for the case when $p$ is a prime power. The question whether the topological Tverberg theorem holds for every p that is not a prime power had been open for long. Very recently, and quite surprisingly, Frick~\cite{Fri} has constructed a counterexample for every non-prime power $p$. His construction is built on work by Mabillard and Wagner~\cite{MaW}. See also \cite{BFZ15} and \cite{AMSW} for further counterexamples.

There is a colourful version of Tverberg theorem. To state it let $n=r(d+1)-1$ and assume that the vertex set $V$ of $\Delta_n$ is partitioned into $d+1$ classes (called colours) and that each colour class contains exactly $r$ vertices. We define $Y_{r,d}$ as the subcomplex of $\Delta_n$ (or $\Delta_n^{(d)}$) consisting of those $\sigma \subset V$ that contain at most one vertex from each colour class. The colourful Tverberg theorem of \v{Z}ivaljevi\'c and Vre\'cica \cite{ZV92} asserts that
$TT(Y_{2p-1,d},d) \geq p$ for prime $p$ which implies that $TT(Y_{4k-1,d},d) \geq k$ for arbitrary $k$. A neat and more recent theorem of Blagojevi\'c, Matschke, and Ziegler ~\cite{BMZ} says that $TT(Y_{r,d},d) = r$ if $r+1$ is a prime, which is clearly best possible. Further information on Tverberg's theorem can be found in Matou\v{s}ek's excellent book \cite{Mato}.

Let $M$ be a matroid (possibly with loops) with rank function $\rho$ on the set $V$. We identify $M$ with the simplicial complex on $V$ whose simplices are the independent sets of $M$. It is well known (see e.g. Theorem 7.8.1 in \cite{B95}) that $M$ is $(\rho(V)-2)$-connected.
Note that both $\Delta_n^{(d)}$ and $Y_{r,d}$ are matroids of rank $d+1$.
In this note we are interested in bounding $TT(M,d)$ for a general matroidal complex
$M$. Let $b(M)$ denote the maximal number of
pairwise disjoint bases in $M$. Our main result is the following
\begin{theorem}
\label{maint} Let $M$ be a matroid of rank $d+1$. Then
$$TT(M,d) \geq \sqrt{b(M)}/4~~.$$
\end{theorem}

In Section \ref{s:conmat} we give a lower bound on the topological
connectivity of the deleted join of matroids. In Section \ref{s:mat}
we use this bound and the approach of \cite{BSS81,ZV92} to prove
Theorem \ref{maint}.

\section{Connectivity of Deleted Joins of Matroids}
\label{s:conmat}

\ \ \ \ We recall some definitions. For a simplicial complex $Y$ on a set $V$ and an element $v \in V$ such that $\{v\} \in Y$, denote the {\it star} and {\it link} of $v$ in $Y$ by
\begin{equation*}
\begin{split}
\st(Y,v)&=\{\sigma \subset V: \{v\} \cup \sigma \in Y\} \\
\lk(Y,v)&=\{\sigma \in \st(Y,v): v \not\in \sigma \}.
\end{split}
\end{equation*}
For a subset $V' \subset V$ let $Y[V']=\{\sigma \subset V': \sigma \in Y\}$ be the induced complex on $V'$.
We regard $\st(Y,v)$, $\lk(Y,v)$ and $Y[V']$ as complexes on the original set $V$ (keeping in mind that not all elements of $V$ have to be vertices of these complexes). Let $f_i(Y)$ denote the number of $i$-simplices in $Y$.
Let $X_1,\ldots,X_k$ be
simplicial complexes on the same set $V$ and let
$V_1,\ldots,V_k$ be $k$ disjoint copies of $V$ with bijections
$\pi_i:V \rightarrow V_i$. The {\it join}  $X_1* \cdots
*X_k$ is the simplicial complex on $ \bigcup_{i=1}^k V_i$
with simplices $\bigcup_{i=1}^k \pi_i(\sigma_i)$ where $\sigma_i \in
X_i$. The {\it deleted join} $(X_1* \cdots
*X_k)_{\Delta}$ is the subcomplex of the join consisting of all simplices $\bigcup_{i=1}^k \pi_i(\sigma_i)$
such that $\sigma_i \cap \sigma_j =\emptyset$ for $1 \leq i \neq j \leq k$. When
all $X_i$ are equal to $X$, we denote their deleted join by
$X^{*k}_{\Delta}$. Note that $\Int_k$ acts freely on $X^{*k}_{\Delta}$ by
cyclic shifts.
\begin{claim}
\label{weakba} Let $M_1,\ldots,M_k$ be matroids on the same
set $V$, with rank functions $\rho_1,\ldots,\rho_k$. Suppose $A_1,\ldots,A_k$ are disjoint subsets of $V$ such
that $A_i$ is a union of at most $m$ independent sets in $M_i$. Then
$Y=(M_1*\cdots *M_k)_{\Delta}$ is $(\lceil\frac{1}{m+1} \sum_{i=1}^k
|A_i|\rceil -2)$-connected.
\end{claim}
\noindent {\bf Proof:} Let
$c=\lceil\frac{1}{m+1} \sum_{i=1}^k |A_i|\rceil-2$.
If $k=1$ then
$\rho_1(V) \geq \left\lceil\frac{|A_1|}{m}\right\rceil$ and hence $Y=M_1$ is $(\left\lceil\frac{|A_1|}{m}\right\rceil-2)$-connected.
For $k \geq 2$ we establish the Claim by induction on $f_0(Y)=\sum_{i=1}^k f_0(M_i)$.
If $f_0(Y)=0$ then all $A_i$'s are empty and the Claim holds. We henceforth assume that $f_0(Y)>0$
and consider two cases:
\\
a) If $M_i=M_i[A_i]$ for all $1 \leq i \leq k$ then $Y=M_1* \cdots *M_k$ is a matroid of rank
$$\sum_{i=1}^k \rho_i(V) \geq \sum_{i=1}^k \left\lceil\frac{|A_i|}{m}\right\rceil \geq \left\lceil\frac{\sum_{i=1}^k |A_i|}{m}\right\rceil.$$
Hence $Y$ is $(\left\lceil\frac{\sum_{i=1}^k |A_i|}{m}\right\rceil-2)$-connected.
\\
b) Otherwise there exists an $1 \leq i_0 \leq k$ such that $M_{i_0} \neq M_{i_0}[A_{i_0}]$. Choose an element
$v \in V-A_{i_0}$ such that $\{v\} \in M_{i_0}$.
Without loss of generality we may assume that $i_0=1$ and that
$v \not\in \bigcup_{i=1}^{k-1} A_i$. Let $S=\bigcup_{i=1}^k V_i$ and let $Y_1=Y[S-\{\pi_1(v)\}]$, $Y_2=\st(Y,\pi_1(v))$.
Then $$Y_1=(M_1[V-\{v\}]*M_2* \cdots *M_k)_{\Delta}.$$
Noting that $f_0(Y_1)=f_0(Y)-1$ and applying the induction hypothesis to the matroids $M_1[V-\{v\}],M_2,\ldots,M_k$ and the sets $A_1,\ldots,A_k$, it follows that $Y_1$ is $c$-connected.
We next consider the connectivity of $Y_1 \cap Y_2$. Write $A_1=\bigcup_{j=1}^t C_j$ where $t \leq m$,  $C_j \in M_1$ for all $1 \leq j \leq t$, and the $C_j$'s are pairwise disjoint.
Since $\{v\} \in M_1$, it follows that there exist $\{C_j'\}_{j=1}^t$ such that $C_j' \subset C_j$, $|C_j'| \geq |C_j|-1$,  and
$C_j' \in \lk(M_1,v)$ for all $1 \leq j \leq t$.
Let
$$
M_i'=
\left\{
\begin{array}{ll}
\lk(M_1,v) & i=1, \\
M_i[V-\{v\}] & 2 \leq i \leq k,
\end{array}
\right.~~
$$
and
\begin{equation*}
A_i'=
\left\{
\begin{array}{ll}
\bigcup_{j=1}^t C_j' & i=1, \\
A_i & 2 \leq i \leq k-1, \\
A_k -\{v\} & i=k.
\end{array}
\right.~~
\end{equation*}
Observe that
$$Y_1 \cap Y_2= \lk(Y,\pi_1(v))=(M_1'* \cdots *M_k')_{\Delta}$$
and that $A_i'$ is a union of at most $m$ independent sets in $M_i'$ for all $1 \leq i \leq k$.
Noting that $f_0(Y_1 \cap Y_2) \leq f_0(Y)-1$ and applying the induction hypothesis to the matroids $M_1',\ldots,M_k'$ and the sets $A_1',\ldots,A_k'$, it follows that $Y_1 \cap Y_2$
is $c'$-connected where
\begin{equation*}
\begin{split}
c'&=\left\lceil\frac{1}{m+1} \sum_{i=1}^k |A_i'|\right\rceil -2 \\
&=
\left\lceil\frac{1}{m+1}\left(\sum_{j=1}^t |C_j'|+
\sum_{i=2}^{k-1}|A_i|+|A_k-\{v\}|\right)\right\rceil-2 \\
&\geq
\left\lceil\frac{1}{m+1}\left(|A_1|-m+
\sum_{i=2}^{k-1}|A_i|+|A_k|-1\right)\right\rceil-2=c-1.
\end{split}
\end{equation*}
As $Y_1$ is $c$-connected, $Y_2$ is contractible and $Y_1 \cap Y_2$ is $(c-1)$-connected, it follows that
$Y=Y_1 \cup Y_2$ is $c$-connected.
{\enp}
\noindent
Let $M$ be a matroid on $V$ with $b(M)=b$ disjoint bases
$B_1,\ldots,B_b$.  Let $I_1 \cup \cdots \cup I_k$ be a partition of
$[b]$ into almost equal parts $\lfloor\frac{b}{k} \rfloor\leq |I_i|
\leq \lceil\frac{b}{k}\rceil$. Applying Claim \ref{weakba} with
$M_1=\cdots=M_k=M$ and $A_i=\cup_{j \in I_i} B_j$, we obtain:
\begin{corollary}
\label{weakb} The connectivity of $M^{*k}_{\Delta}$ is at least
$$\frac{b\rho(V)}{\lceil\frac{b}{k}\rceil+1}-2~~.$$
\end{corollary}
\noindent
We suggest the following:
\begin{conjecture}
\label{deljoin} For any $k \geq 1$ there exists an $f(k)$ such that
if $b(M) \geq f(k)$ then $M^{*k}_{\Delta}$ is  $(k
\rho(V)-2)$-connected.
\end{conjecture}
\noindent {\bf Remark:} Let $M$ be the rank $1$ matroid on $m$ points $M=\Delta_{m-1}^{(0)}$. The chessboard complex $C(k,m)$ is the $k$-fold deleted join $M^{*k}_{\Delta}$.
Chessboard complexes play a key role in the works of \v{Z}ivaljevi\'c and Vre\'cica \cite{ZV92} and
Blagojevi\'c, Matschke, and Ziegler ~\cite{BMZ} on the colourful Tverberg theorem.
Let $k \geq 2$. Garst \cite{Garst} and \v{Z}ivaljevi\'c and Vre\'cica \cite{ZV92}
proved that $C(k,2k-1)$ is $(k-2)$-connected.
On the other hand, Friedman and Hanlon \cite{FH98} showed that $\tilde{H}_{k-2}(C(k,2k-2);\Rat) \neq 0$,
so $C(k,2k-2)$ is not $(k-2)$-connected. This implies that the function $f(k)$ in Conjecture \ref{deljoin} must satisfy $f(k) \geq 2k-1$.

\section{A Tverberg Type Theorem for Matroids}\label{s:mat}

 We recall some well-known topological facts
(see \cite{BSS81}). For $m \geq 1, k \geq 2$ we identify the sphere
$S^{m(k-1)-1}$ with the space
$$\left\{(y_1,\ldots,y_k) \in (\Rea^m)^k: \sum_{i=1}^k
|y_i|^2=1~,~\sum_{i=1}^k y_i=0 \in \Rea^m ~\right\}~~.$$ The cyclic shift on this
space defines a $\Int_k$ action on $S^{m(k-1)-1}$. The action is
free for prime $k$.

The {\it $k$-fold deleted product} of a space $X$ is the
$\Int_k$-space given by
$$X^k_D=X^k-\{(x,\ldots,x)
\in X^k : x \in X \}~.$$ For $m \geq 1$ define a $\Int_k$-map
$$\phi_{m,k}: (\Rea ^m)^k_D \rightarrow S^{m(k-1)-1}$$
by $$\phi_{m,k}(x_1,\ldots,x_k)=\frac{(x_1-\frac{1}{k}\sum_{i=1}^k
x_i,\ldots,x_k-\frac{1}{k}\sum_{i=1}^k x_i)}{(\sum_{j=1}^k |
x_j-\frac{1}{k}\sum_{i=1}^k x_i|^2)^{1/2}}~.$$ We'll also need the
following result of Dold \cite{Dold} (see also Theorem 6.2.6 in \cite{Mat}):
\begin{theorem}[Dold]
\label{dold} Let $p$ be a prime and suppose $X$ and $Y$ are free
$\Int_p$-spaces such that $\dim Y=k$ and $X$ is $k$-connected. Then
there does not exist a $\Int_p$-map from $X$ to $Y$.
\end{theorem}

\noindent {\bf Proof of Theorem \ref{maint}:} Let $M$ be a matroid
on the vertex set $V$, and let $f:M \rightarrow \Rea^d$ be a
continuous map. Let $b=b(M)$ and choose a prime $\sqrt{b}/4 \leq p
\leq \sqrt{b}/2$. We'll show that there exist disjoint simplices
(i.e. independent sets) $\sigma_1,\ldots,\sigma_p \in M$ such that
$\bigcap_{i=1}^p f(\sigma_i) \neq \emptyset$. Suppose for contradiction
that $\bigcap_{i=1}^p f(\sigma_i)= \emptyset$ for all such choices of
$\sigma_i$'s. Then $f$ induces a continuous $\Int_p$-map
$$f_*:M^{*p}_{\Delta} \rightarrow (\Rea^{d+1})^p_D$$
as follows. If $x_1,\ldots,x_p$ have pairwise disjoint supports in
$M$ and $(t_1,\ldots,t_p) \in \Rea^p_+$ satisfies $\sum_{i=1}^p
t_i=1$ then
$$f_*(t_1\pi_1(x_1)+ \cdots +t_p \pi_p(x_p))=
(t_1,t_1f(x_1),\ldots,t_p,t_pf(x_p)) \in (\Rea^{d+1})^p_D~~.$$ Hence
$\phi_{d+1,p}f_*$ is a $\Int_p$-map between the free $\Int_p$-spaces
$M^{*p}_{\Delta}$ and $S^{(d+1)(p-1)-1}$. This however contradicts
Dold's Theorem since by Corollary \ref{weakb} the connectivity of
$M^{*p}_{\Delta}$ is at least
$$\frac{b(d+1)}{\lceil\frac{b}{p}\rceil+1} -2 \geq (d+1)(p-1)-1~$$
by the choice of $p$.
{\enp}
\ \\ \\
{\bf Acknowledgements.} Research of Imre B\'ar\'any was partially supported by ERC advanced grant 267165, and by Hungarian National grant K 83767. Research of Gil Kalai was
supported by ERC advanced grant 320924. Research of Roy Meshulam is supported by ISF and GIF grants.

\vskip.9cm
\noindent Authors' addresses:\\[2mm]
Imre B\'{a}r\'{a}ny\\
R\'{e}nyi Institute, Hungarian Academy of Sciences\\
POB 127, 1364 Budapest, Hungary\\
and\\
Department of Mathematics, University College London\\
Gower Street, London, WC1E 6BT, UK\\
E-mail: barany@renyi.hu\\[3mm]
Gil Kalai\\
Einstein Institute of Mathematics, Hebrew University\\
Jerusalem 9190, Israel\\
E-mail: kalai@math.huji.ac.il \\[3mm]
Roy Meshulam\\
Department of Mathematics,  Technion\\
Haifa 32000,  Israel\\
E-mail: meshulam@math.technion.ac.il

\end{document}